\documentclass[12pt]{article}
\usepackage[table]{xcolor}
\usepackage{graphicx}
\usepackage{amsmath,amsfonts,amssymb}
\usepackage{graphicx}
\usepackage{multirow}
\usepackage{tikz,pgf}
\usepackage{url}
\usepackage{amsmath}
\usepackage{graphicx}
\usepackage{epsfig}
\usepackage{epstopdf}
\usepackage{color}
\usepackage{enumitem}
\def\tto{\;{\lower 1pt \hbox{$\rightarrow$}}\kern -10pt
\hbox{\raise 2pt \hbox{$\rightarrow$}}\;}

\def\ra{\rangle}
\def\la{\langle}

\def\epsilon{\varepsilon}

\def\h{\hfill\Box}
\def\R{\mathbb{R}}

\def\N{\mathbb{N}}

\def\gph{\mbox{\rm gph}\,}

\def\epi{\mbox{\rm epi}\,}

\def\dom{\mbox{\rm dom}\,}

\def\h{\hfill\square}

\setlist[enumerate,1]{itemsep=0.0ex,parsep=0.5ex,label={\rm(\alph*)},leftmargin=*,align=left}
\newcounter{lk}

\newcommand{\conv}{{\rm conv}}

\usepackage[colorlinks=true]{hyperref}

\begin{document}
\begin{center}
{\sc\bf NASH EQUILIBRIUM AND MINIMAX THEOREMS\\ VIA VARIATIONAL TOOLS OF CONVEX ANALYSIS}\\[3ex]
{\sc Nguyen Xuan Duy Bao\footnote{Department of Applied Mathematics, University of Technology, Ho Chi Minh City, Vietnam \\ Vietnam National University, Ho Chi Minh City, Vietnam  (nxdbao@hcmut.edu.vn)}, Boris S. Mordukhovich\footnote{Department of Mathematics, Wayne State University, Detroit, Michigan 48202, USA (boris@math.wayne.edu). Research of this author was partly supported by the USA National Science Foundation under grant DMS-1808978, by the Australian Research Council under grant DP-190100555, and by Project 111 of China under grant D21024.}, Nguyen Mau Nam\footnote{Fariborz Maseeh Department of Mathematics and Statistics, Portland State University, Portland, OR
97207, USA (mnn3@pdx.edu).}}
\end{center}
\small{\bf Abstract.} In this paper, we first provide a simple variational proof of the existence of Nash equilibrium in Hilbert spaces by using optimality conditions in convex minimization and Schauder's fixed-point theorem. Then applications of convex analysis and generalized differentiation are given to the existence of Nash equilibrium and extended versions of von Neumann's minimax theorem in locally convex topological vector spaces. Our analysis in this part combines generalized differential tools of convex analysis with elements of fixed point theory revolving around Kakutani's fixed-point theorem and related issues.\\[1ex]
\normalsize
\begin{center}
{\sc\bf Dedicated to Juan Enrique Mart\'inez-Legaz\\ in friendship and great respect}
\end{center}

\newtheorem{Theorem}{Theorem}[section]
\newtheorem{Proposition}[Theorem]{Proposition}
\newtheorem{Remark}[Theorem]{Remark}
\newtheorem{Lemma}[Theorem]{Lemma}
\newtheorem{Corollary}[Theorem]{Corollary}
\newtheorem{Definition}[Theorem]{Definition}
\newtheorem{Example}[Theorem]{Example}
\renewcommand{\theequation}{\thesection.\arabic{equation}}
\normalsize

\section{Introduction and Overview}\label{intro}

Game theory has its origin in the work of von Neumann on the {\em minimax theorem} \cite{Nm1928}. This theorem establishes sufficient conditions under which the inequality 
\begin{equation}\label{eq1.1}
\sup_{x\in \Omega_1}\inf_{y\in \Omega_2}f(x,y)\leq \inf_{y\in \Omega_2}\sup_{x\in \Omega_1}f(x,y)
\end{equation}
holds as an equality, where $f\colon \Omega_1\times \Omega_2\to \mathbb R$ is a function defined on the Cartesian product of two nonempty  sets $\Omega_1$ and $\Omega_2$.

In the context of two-person noncooperative zero-sum games, $\Omega_i$ is the set of strategies of the $i$th player for $i=1,2$, and $f(x, y)$ represents the  payoff to the first player choosing strategy $x$ and the second player choosing strategy $y$. Because this is a zero-sum game, the payoff to the second player is $-f(x, y)$. Given strategy $x\in\Omega_1$ of the first player, the second player uses strategy $y\in \Omega_2$ to minimize the first player's payoff. The expression $\inf_{y\in \Omega_2}f(x, y)$ represents  the  minimum payoff of the first player for playing strategy $x$. Hence the expression on the left-hand side of \eqref{eq1.1} signifies the {\em best worst-case scenario} for the first player who seeks strategy $x$ that maximizes this minimum payoff. Similarly, $\sup_{x\in \Omega_1}f(x, y)$ represents the maximum payoff of the first player for each strategy $y\in \Omega_2$ chosen by the second player. Therefore, the expression on the right-hand side of \eqref{eq1.1} signifies the {\em worst maximum payoff} for the first player when the second player seeks strategy $y\in \Omega_2$ to minimize this maximum payoff. From the perspective of the first player, inequality \eqref{eq1.1} means that {\em best worst} is always less than or equal to {\em worst best}.

It turns out that if there exists an {\em saddle point} $(\bar x,\bar y)\in \Omega_1\times \Omega_2$ for $f$ such that
\begin{equation*}
f(x, \bar y)\leq f(\bar x,\bar y)\leq f(\bar x, y)\; \mbox{\rm whenever }\;x\in \Omega_1\;\mbox{ and }\;y\in \Omega_2,
\end{equation*}
then \eqref{eq1.1} holds as an equality. At the saddle point $(\bar x_,\bar y)$, if the first player decides to change strategy from $\bar x$ to $x$, then $f(x,\bar y)\leq f(\bar x,\bar y)$ means that this player cannot improve their payoff.  Meanwhile, if the second player chooses to change strategy from $\bar y$ to $y$, then $f(\bar x,\bar y)\leq f(\bar x, y)$, and so $-f(\bar x, y)\leq -f(\bar x,\bar y)$. Thus the second player cannot improve their payoff as well. Therefore, {\em neither player can unilaterally change their strategy to improve their payoff} at the equilibrium point. It is proved in his landmark paper by von Neumann \cite{Nm1928} that every two-person noncooperative  zero-sum game in which the payoff function $f(x, y)=x^TAy$, where $A$ is an $n\times m$ matrix, admits a saddle point in  {\em mixed strategies} represented by two nonempty compact convex sets $\Omega_1\subset \mathbb R^n$ and $\Omega_2\subset\mathbb R^m$. This point is also referred to as an {\em equilibrium} or {\em saddle point} for $f$ over $\Omega_1\times\Omega_2$.

A natural question arises on how to extend von Neumann's work on the minimax theorem to the case of two-person noncooperative games.  This question was first addressed by Nash in \cite{nashthesis} who introduced and investigated the notion of {\em Nash equilibrium} that has been broadly applied to economics, ecology, various fields of science and technology, etc. Among a great many important developments in these directions, we mentioned the remarkable paper by Bonenti and Mart\'inez-Legaz \cite{je} who achieved a complete characterization of the existence of saddle points in a fairly general setting with applications to convex programming. 

Now we recall the fundamental  notion of Nash equilibrium. Consider the two functions $u_1, u_2\colon \Omega_1\times \Omega_2\to \mathbb R$ defined on the Cartesian product of nonempty sets $\Omega_1$ and $\Omega_2$, where $u_i$ indicates the payoff function, while $\Omega_i$ represents the set of strategies of the $i$th player as $i=1, 2$. Then a {\em Nash equilibrium point} for $(u_1,u_2)$ is a pair $(\bar x,\bar y)\in \Omega_1\times \Omega_2$ satisfying the following conditions:
\begin{equation}\label{nash}
\begin{array}{ll}
u_{1}(x,\bar y) \leq  u_{1}(\bar x,\bar y) \ \;  \text{for all} ~ x \in \Omega_{1}, \\
u_{2}(\bar x, y) \leq u_{2}(\bar x,\bar y) \ \; \text{for all} ~ y \in \Omega_{2}.
\end{array}
\end{equation}

In the original work by Nash \cite{nashthesis}, the existence of an equilibrium point is proved when $u_1:=x^TAy$ and $u_2:=x^TBy$, with $A$ and $B$ standing $n\times m$ matrices, in the class of mixed strategies belonging to nonempty compact convex sets $\Omega_1\subset \mathbb R^n$ and $\Omega_2\subset\mathbb R^m$. Similarly to the case of saddle points, no players can benefit at the equilibrium point $(\bar x,\bar y)$ by changing their strategies while the other players' strategies remain unchanged. Observe to this end that the notion of Nash equilibrium is a generalization of the saddle point notion because when $u_2=-u_1$, any Nash equilibrium point $(\bar x,\bar y)\in \Omega_1\times \Omega_2$ for $u_1$ and $u_2$  is also a saddle point for $u_1$.

{\em Convexity} plays a crucial role in both von Neumann's and Nash's works, while their proofs of the existence of equilibrium points use rather involved arguments. An effort to apply modern theory of convex analysis with shedding new lights on Nash equilibrium was accomplished by Rockafellar in \cite{Rockafellar2012} who used Brouwer's fixed-point theorem together with elements of convex analysis and theory of maximal monotone operators to prove the existence of Nash equilibrium when the functions $u_i$ involved are differentiable and $\Omega_i$ for $i=1, 2$ are nonempty compact convex sets in finite dimensions. In Section~3 below, we provide a simple {\em variational} proof of the existence of Nash equilibrium, not merely in finite-dimensional but in infinite-dimensional Hilbert spaces, by using optimality conditions in convex programming and Schauder's fixed-point theorem. Subsequently, in Section 4, Kakutani's fixed-point theorem for multifunctions is employed to establish a novel result on the existence of Nash equilibrium in cases where $u_i$ may not be differentiable. This extension enables us to develop a new proof from convex analysis for the existence of von Neumann's saddle points in the context of two-player noncooperative zero-sum games with nondifferentiable payoff functions. In Section 5, we revisit variations of the minimax theorem in locally convex topological vector (LCTV) spaces and offer direct proofs based on the convex separation theorem. For the reader's convenience, in Section~2 we overview the needed preliminaries from convex analysis and fixed point theory.
 
\section{Preliminaries from  Convex Analysis and Fixed Point Theory in LCTV Spaces}

We start with some elements of convex analysis used below, while referring the reader to the book  \cite{nambook22, zalinescu02} for more details in LCTV spaces. To those who are more interested in finite-dimensional theory, we recommend the books \cite{nambook23, rockbook}. Having a (Hausdorff) LCTV space $X$, denote by $X^*$ its topological dual and use  $\la x^*, x\ra$ to indicate canonical pairing between $x^*\in X^*$ and $x\in X$. The notation $\mathbb{R}^n_+$ signifies the nonnegative orthant in $\mathbb{R}^n$. From now on, we assume that all vector spaces under consideration are over the real line $\R$.

Recall that a subset $\Omega$ of an LCTV space $X$ is {\em convex} if, along with any two points $a, b\in\Omega$, it contains the line segment defined by
$$
[a, b]: = \big\{\lambda a + (1-\lambda) b \; \big |\;  0\leq \lambda\leq 1\big\}.
$$
In other words, $\Omega$ is convex if and only if $\lambda a + (1-\lambda)b \in \Omega$ for all $a, b \in \Omega$ and $\lambda \in [0,1]$. The {\em convex hull}  $\conv(\Omega)$ of $\Omega$ is the collection of all convex combinations of points in $\Omega$. The {\em indicator function}  $\delta(\cdot; \Omega)\colon X\to (-\infty, \infty]$ associated with $\Omega$ is 
$$
\delta(x; \Omega): = \begin{cases} 
0 & \text{if } x \in \Omega, \\
\infty & \text{otherwise}.
\end{cases}
$$
Given an element $\bar x\in \Omega$, the {\em normal cone} to $\Omega$ at $\bar x$ is
\begin{equation*}
N(\bar x;\Omega): = \big\{x^*\in  X^*\; \big |\; \langle x^*, x-\bar x\rangle \leq  0 \ \ \; \mbox{\rm for all }\;x \in \Omega\big\}.
\end{equation*}

Turning now to extended-real-valued functions $f\colon X\rightarrow (-\infty, \infty]$, recall that $f$ is \emph{lower semicontinuous} on a subset $\Omega$ of $X$ if for each $\alpha \in \mathbb{R}$, the lower level set $L_{\alpha}(f): = \{x \in \Omega \mid f(x) \leq \alpha\}$ is  closed in $\Omega$ with the induced topology. Symmetrically, $f\colon X\to [-\infty, \infty)$ is \emph{upper semicontinuous} if $-f$ is lower semicontinuous on $\Omega$. The {\em continuity} of $f$ on a subset $\Omega$ of $X$ means that $f(x)\in \R$ for every $x\in \Omega$, and for any $x_0\in \Omega$ and any net $\{x_\tau\}\subset \Omega$ converging to $x_0$ we have $f(x_\tau)\to f(x_0)$.

Next we formulate the notions of (generalized) convexity and concavity for extended-real-valued functions defined on convex sets $\Omega\subset X$ used in what follows. A function $f\colon \Omega \rightarrow (-\infty, \infty]$ is {\em convex} if
$$
f\big(\lambda x + (1-\lambda)y\big)\le\lambda f(x) + (1-\lambda) f(y)\; \ \mbox{\rm for all }\;x, y \in \Omega\; \ \mbox{\rm and }\;0<\lambda<1.$$
We say that $f\colon \Omega \rightarrow (-\infty, \infty]$ is {\em quasiconvex} if
$$
f\big(\lambda x + (1-\lambda)y\big) \leq \max\big\{f(x), f(y)\big\} \; \ \mbox{\rm for all }\;x, y \in \Omega\; \ \mbox{\rm and }\;0<\lambda<1.
$$
Symmetrically, $f\colon \Omega\to[-\infty, \infty)$ is called {\em concave} (resp.\ {\em quasiconcave}) if the function $-f$ is convex (resp.\ quasiconvex).  

It is well known that $f$ is convex if and only if its epigraphical set  $\epi f:=\{(x,\alpha)\in X\times\R \mid\alpha\ge f(x)\}$ is convex, and it is quasiconvex if and only if each lower level set $L_\alpha(f)$ is convex. Furthermore,  any real-valued function $f\colon\R^n\to\R$, which is either convex or concave, is locally Lipschitz continuous.

Now we define a central generalized differential notion of convex analysis. Given a convex function $f\colon X\to (-\infty, \infty]$, define the {\em effective domain} of $f$ by
\begin{equation*}
    \dom f:=\big\{x\in X\; \big |\; f(x)<\infty\big\}.
\end{equation*}
For a point $\bar x\in \dom f$, a dual element $x^*\in X^*$ is a {\em subgradient} of $f$ at $\bar x$ if
\begin{equation*}
\langle x^*, x-\bar x\rangle \leq f(x)-f(\bar x)\; \mbox{\rm for all }\;x \in X.
\end{equation*}
The collection of all subgradients of $f$ at $\bar x$ is called the {\em subdifferential} of $f$ at this point and is denoted by $\partial f(\bar x)$. We set $\partial f(\bar x):=\emptyset$ if $f(\bar x)=\infty$. Let $f\colon X\times Y\to (-\infty, \infty]$ be a function of two variables $(x, y) \in X \times Y$ which is convex with respect to $x$. For a point $(\bar x,\bar y)\in X \times Y$, the {\em partial subdifferential} of $f$ with respect to $x$ at this point denoted by $\partial_x f(\bar x,\bar y)$ is defined as the subdifferential of the function $f(\cdot, \bar y)$ at $\bar x$.

Given further a cost/objective function $f\colon X \rightarrow \mathbb{R}$ and a constraint set $\Omega\subset X$, consider the following optimization problem $(\mathrm{P})$ :
\begin{center}
minimize $f(x)$\ \; subject to ~ $x \in \Omega$
\end{center}
for which $\bar x\in \Omega$ is an {\em optimal solution} (global minimizer) if $f(x_0)\le f(x)$ whenever  $x \in \Omega$. Here is a 
well-known characterization of optimal solutions to convex problems of constrained minimization.

\begin{Proposition}[\bf characterizing optimal solutions to convex problems]\label{prop2.4}
Let $f\colon X \rightarrow \R$ be a continuous convex function on an LCTV space $X$, and let $\Omega$ be a convex subset of $X$ with $\bar x\in \Omega$. Then $\bar x$ is an optimal solution to $(\mathrm{P})$ if and only if
$$
0 \in \partial f(\bar x) + N(\bar x;\Omega).
$$
In the case where $f$ is $($Fr\'echet$)$ differentiable at $\bar x$, the inclusion above can be written in the form $-\nabla f(\bar x) \in N(\bar x;\Omega)$.
\end{Proposition}

In the case where $X$ is a normed space, the {\em distance function} associated with a nonempty set $\Omega\subset X$ is given by
\begin{equation} \label{eq2.1}
d(x;\Omega): = \inf \big\{\|x-\omega\|\; \big |\;  \omega \in \Omega\big\}, \ \; x \in X.
\end{equation}	
For each element $x \in X$, the {\em metric projection} from $x$ to $\Omega$ is defined via \eqref{eq2.1} as
\begin{equation} \label{eq2.2}
\mathcal{P}(x;\Omega): = \big\{\omega \in \Omega \; \big |\; \|x-\omega\| = d(x;\Omega)\big\}.
\end{equation}
A remarkable consequence of convexity in Hilbert spaces is the following.

\begin{Proposition}[\bf projections to convex sets]\label{hilb-proj}
Let $\Omega$ be a nonempty closed convex subset of a Hilbert space $H$. Then for any $x \in H$, the metric projection $\mathcal{P}(x;\Omega)$ is a singleton. Furthermore,  $\bar\omega\in \mathcal{P}(x;\Omega)$ if and only if $\bar\omega\in \Omega$ and
$$
\langle x-\bar\omega, \omega-\bar\omega\rangle \leq 0 ~ \text{for all } ~ \omega \in \Omega.
$$
\end{Proposition}

Therefore, the assumptions of Proposition~\ref{hilb-proj} ensure that the projection mapping associated with the set $\Omega$ is single-valued. The next proposition asserts that this mapping is Lipschitz continuous with constant $\ell=1$ on $H$.

\begin{Proposition}[\bf Lipschitz continuity of projections]\label{prop2.6} Let $\Omega\subset H$ satisfy the assumptions of Proposition~{\rm\ref{hilb-proj}}.  Then the projection mapping $\mathcal{P}(\cdot;\Omega)\colon H\to H$ is Lipschitz continuous on $H$ with the global Lipschitz constant $\ell=1$, i.e.,
\begin{equation*}
\left\|\mathcal{P}(x_{1};\Omega) - \mathcal{P}(x_{2};\Omega)\right\| \leq \left\|x_{1}-x_{2}\right\| \ \;  \mbox{\rm for all }\; x_{1}, x_{2} \in H.
\end{equation*}
\end{Proposition}

Yet another useful property relates projections and normals to convex sets.

\begin{Proposition}[normals via projections]\label{prop2.7}
Let $\Omega\subset H$ satisfy the assumptions of  Proposition~{\rm\ref{hilb-proj}} , and let $\bar x\in \Omega$. Then $v \in N(\bar x;\Omega)$ if and only if $\bar x= \mathcal{P}(\bar x+v;\Omega)$.
\end{Proposition}

Next we formulate two fundamental {\em fixed-point} theorems. The first one is Schauder's results for continuous single-valued mappings in LCTV spaces, which is an infinite-dimensional extension of the classical Brouwer fixed-point theorem established in finite dimensions; see \cite{zeidler} for more details and discussions.

\begin{Theorem}[\bf Schauder's fixed-point theorem]\label{theo2.8}
Let $X$ be an LCTV space, and let $f\colon \Omega\to \Omega$ be a continuous mapping from a nonempty compact convex set $\Omega$ to itself. Then $f$ has a fixed point, i.e., there exists $\bar x\in \Omega$ such that $f(\bar x)=\bar x$.
\end{Theorem}

We'll also need fixed-point results for set-valued mappings/multifunctions $F\colon\Omega\tto\Omega$, i.e., mappings taking values in the collections of subsets of $\Omega$. Here is an extension of Theorem~\ref{theo2.8}, which goes back to Kakutani; see \cite{sharpirobook}.

\begin{Theorem}[Kakutani's fixed-point theorem]\label{theo2.9}
Let $X$ be a LCTV space, and let $\Omega$ be a nonempty compact convex subset of $X$. Suppose that $F\colon \Omega\tto \Omega$ is a multifunction satisfying the properties:

{\bf(a)} The sets $F(x)$ is  nonempty and convex for all $x\in \Omega$.

{\bf(b)} The graph $\gph F:=\{(x, y)\; |\; x\in \Omega,\;  y\in F(x)\}$ is closed in $X\times X$.
     
Then $F$ has a fixed point, i.e., there exists $\bar x\in \Omega$ such that $\bar x\in F(\bar x)$.
\end{Theorem}

\section{Variational Approach to Nash Equilibrium}

This section offers a brief introduction to noncooperative game theory with two instructive illustrative examples. Then we give a simple variational proof of the classical existence theorem for Nash equilibrium in finite-dimensional spaces and establish its extension to Hilbert spaces by using optimality conditions in convex minimization and  Schauder's fixed-point theorem.

\begin{Example}[\bf Nash equilibrium with pure strategies]\label{ex3.1} {\rm 
In a scenario involving two competing tech companies, Company X and Company Y, both companies must decide whether to invest in research and development (R\&D). If both companies decide to allocate resources towards R\&D, then each company will experience a positive outcome of 50 units. If only one company chooses to invest in R\&D, it will gain 100 units while the other company will suffer a loss of 50 units as a consequence of a decrease in market shares. However, if neither company opts to invest in R\&D, then both companies will receive a net benefit of 0 units. The payoff function of each company is represented in the payoff matrix below:

\begin{center}
\begin{tabular}{cc|lr|lr|}
&\multicolumn{1}{c}{} & \multicolumn{4}{c}{Company Y} \\
& \multicolumn{1}{c}{} & \multicolumn{2}{c}{\hspace{0.25cm}Do} & \multicolumn{2}{c}{\hspace{0.25cm}Don't} \\\cline{3-6}
\multirow{4}*{\rotatebox{90}{Company X}}  & \multirow{2}*{Do} & & 50 & & -50 \\
& & 50 & & 100 & \\\cline{3-6}
& \multirow{2}*{Don't} & & 100 & & 0 \\
& & -50 & & 0 & \\\cline{3-6}
\end{tabular}
\end{center}
In this example, the Nash equilibrium occurs when both companies implement their {\em pure strategies} by choosing to invest in R\&D. This is because the choice to invest yields a better or equal payoff regardless of the opponent's action:
\begin{itemize}
\item If Company X invests, it either gains 50 units (if Company Y also invests) or 100 units (if Company Y does not invest). Therefore, investing is better compared to not investing (where it would gain 0 units if Company Y also does not invest, or lose 50 units if Company Y invests).
\item The same rationale applies to Company Y.
\end{itemize}
Thus both investing in R\&D are Nash equilibria since both companies experience deceasing payoff by unilaterally changing their strategy from investing to not investing.}
\end{Example}

\begin{Example}[Nash equilibrium with mixed strategies] \label{ex3.2} {\rm 
Let us consider another simple two-person game called \emph{matching pennies}. Suppose that Player I and Player II each has a penny. Each player must secretly turn the penny to heads or tails and then reveal their choices simultaneously. Player I wins the game and gets Player II's penny if both coins show the same face (heads-heads or tails-tails). Conversely, if  the coins show different faces (heads-tails or tails-heads), Player II wins the game and gets Player I's penny. The payoff function of each player is represented in the payoff matrix given below.
\begin{center}
\begin{tabular}{cc|lr|lr|}
   & \multicolumn{1}{c}{} & \multicolumn{4}{c}{Player II} \\
   & \multicolumn{1}{c}{} & \multicolumn{2}{c}{head} & \multicolumn{2}{c}{tail} \\\cline{3-6}
   \multirow{4}*{\rotatebox{90}{Player I}} & \multirow{2}*{head} & & $-1$ & & 1 \\
   & & 1 & & $-1$ & \\\cline{3-6}
   & \multirow{2}*{tail} & & 1 & & $-1$ \\
   & & $-1$ & & 1 & \\\cline{3-6}
\end{tabular}
\end{center}
In this game it is not hard to see that no pure strategy Nash equilibrium exists. Suppose that Player I is playing with a mind reader who knows Player I's choice of faces. If Player I decides to turn heads, then Player II knows about it and chooses turn tails to win the game. In the case where Player I decides to turn tails, then Player II chooses to turn heads and again wins the game. Thus to have a fair game, each player decides to randomize their strategy by, for instance, tossing the coin instead of putting the coin down. We describe the new game as follows.
\begin{center}
\begin{tabular}{cc|lr|lr|}
   & \multicolumn{1}{c}{} & \multicolumn{4}{c}{Player II} \\
   & \multicolumn{1}{c}{} & \multicolumn{2}{c}{head, $q_1$} & \multicolumn{2}{c}{tail, $q_2$} \\\cline{3-6}
   \multirow{4}*{\rotatebox{90}{Player I}} & \multirow{2}*{head, $p_1$} & & $-1$ & & 1 \\
   & & 1 & & $-1$ & \\\cline{3-6}
   & \multirow{2}*{tail, $p_2$} & & 1 & & $-1$ \\
   & & $-1$ & & 1 & \\\cline{3-6}
\end{tabular}
\end{center}
In this new game, Player I uses a coin randomly with probability of coming up heads $p_{1}$ and probability of coming tails $p_{2}$, where $p_{1}+p_{2}=1$. Similarly, Player II uses another coin randomly with probability of coming up heads $q_{1}$ and probability of coming tails $q_{2}$, where $q_{1}+q_{2}=1$. The new strategies are now called mixed strategies while the original ones are called pure strategies. Denote by $A$ the matrix that represents the payoff function of Player I, and by $B$ the matrix that represents the payoff function of Player II:
	$$
	A ~ = ~ \begin{array}{|c|c|}
		\hline	1 & -1 \\
		\hline	-1 & 1 \\
		\hline
	\end{array} 
	\quad \text{and} \quad
	B ~ = ~ \begin{array}{|c|c|}
		\hline	-1 & 1 \\
		\hline 	1 & -1 \\
		\hline
	\end{array}
	$$
Suppose that Player II uses mixed strategy $\{q_{1}, q_{2}\}$, then Player I's expected payoff for playing heads is
	$$u_{H}(q_{1}, q_{2}) = q_{1}-q_{2} = 2q_{1}-1.$$
Similarly, Player I's expected payoff for playing tails is
	$$u_{T}(q_{1}, q_{2}) = -q_{1}+q_{2} = -2 q_{1}+1.$$
Thus Player I's expected payoff for playing the mixed strategy $\{p_{1}, p_{2}\}$ is
	$$u_{1}(p,q) = p_{1}u_{H}(q_{1}, q_{2}) + p_{2} u_{T}(q_{1}, q_{2}) = p_{1}(q_{1}-q_{2}) + p_{2}(-q_{1}+q_{2}) = p^{T}Aq$$
where $p = [p_{1}, p_{2}]^{T}$ and $q = [q_{1}, q_{2}]$.

By the same arguments, if Player I chooses mixed strategy $\{p_{1}, p_{2}\}$, then Player II's expected payoff for playing mixed strategy $\{q_{1}, q_{2}\}$ is
	$$u_{2}(p, q) = q^{T}Bp.$$
In this new game, an element $(\bar p,\bar q)$ is a Nash equilibrium if
\begin{align*}
&u_{1}(p,\bar q) \leq u_{1}(\bar p,\bar q) \ \; \text{for all} ~ p \in \Delta, \\
&u_{2}(\bar p, q) \leq u_{2}(\bar p,\bar q) \ \; \text{for all} ~ q \in \Delta,
\end{align*}
where $\Delta: = \left\{(p_{1}, p_{2}) \in \mathbb{R}^2 \mid p_{1} \geq 0, p_{2} \geq 0, p_{1}+p_{2}=1 \right\}$ is a nonempty compact convex subset of $\mathbb{R}^{2}$}.
\end{Example}

In \cite{nashthesis}, Nash proved the existence of an equilibrium point \eqref{nash} in the class of mixed strategies in finite-dimensional spaces covering, in particular, the game in Example~\ref{ex3.2}. His proof was based on Brouwer's fixed-point theorem as well as rather involved arguments from polyhedrality theory. In \cite{Rockafellar2012}, Rockafellar suggested another approach to prove the classical Nash equilibrium theorem in finite dimensions by using, along with Brouwer's fixed-point theorem, theory of maximal monotone operators. We now develop a new, much simplified proof---first for the classical Nash theorem, and then for its extension to Hilbert spaces---based on characterizing minimizers in standard problems of convex constrained optimization with smooth objectives.

\begin{Theorem}[classical Nash equilibrium theorem] \label{theo3.3}
Consider a two-person game $\{\Omega_{i}, u_{i}\}$, $i=1,2$, where $\Omega_{1} \subset \mathbb{R}^{m}$ and $\Omega_{2} \subset \mathbb{R}^{n}$ are nonempty compact convex sets, and where the payoff functions $u_{i}\colon \Omega_{1} \times \Omega_{2} \rightarrow \mathbb{R}$ are given by
$$
u_{1}(p, q): = p^{T}Aq\;\mbox{ and }\;u_{2}(p, q): = q^{T}Bp
$$
with $A$ and $B$ standing for $m \times n$ and $n \times m$ matrices, respectively. Then this game admits a Nash equilibrium.
\end{Theorem}
{\bf Proof.} It follows from definition \eqref{nash} that a pair $(\bar p,\bar q) \in \Omega_{1} \times \Omega_{2}$ is a Nash equilibrium of the game under consideration if and only if
\begin{equation}\label{nash1}
\begin{array}{ll}
&-u_{1}(p,\bar q) \geq -u_{1}(\bar p,\bar q) \ \; \text{for all} ~ p \in \Omega_{1}, \\
&-u_{2}(\bar p, q) \geq -u_{2}(\bar p,\bar q) \ \; \text{for all} ~ q \in \Omega_{2}.
\end{array}
\end{equation}
Proposition~\ref{prop2.4} tells us that the above inequalities hold if and only if we have the normal cone inclusions
\begin{equation} \label{eq3.1}	
\nabla_{p} u_{1}(\bar p,\bar q) \in N(\bar p;\Omega_{1}) \ \; \text{and} \ \; \nabla_{q} u_{2}(\bar p,\bar q) \in N(\bar q;\Omega_{2}).
\end{equation}
By using the structures of $u_{i}$, conditions \eqref{eq3.1} can be equivalently expressed as
$$
A\bar q\in N(\bar p;\Omega_{1}) \ \; \text{and} \ \;  B\bar p\in N(\bar q;\Omega_{2}),
$$
which can be rewritten in the form
$$
(A\bar q, B\bar p)\in N(\bar p;\Omega_1)\times N(\bar q;\Omega_2)=N\big((\bar p,\bar q);\Omega_1\times\Omega_2\big).
$$
Employing Proposition~\ref{prop2.7} on the metric projection, the latter is equivalent to
\begin{equation} \label{eq3.2}
(\bar p,\bar q) = \mathcal{P}\big((\bar p,\bar q) + (A\bar q, B\bar p);\Omega\big)\;\mbox{ with }\;\Omega:=\Omega_1\times\Omega_2.
\end{equation}
Defining now the mapping $\Phi\colon \Omega \rightarrow \Omega$ by
$$
\Phi(p, q): = \mathcal{P}\big((p, q)+(Aq, Bp);\Omega\big)
$$
and implementing another elementary projection property from Proposition~\ref{prop2.6} ensure that the mapping $\Phi$ is (Lipschitz) continuous while the set $\Omega$ is nonempty, convex, and compact. According to the classical fixed-point theorem by Brouwer in finite dimensions, the mapping $\Phi$ has a fixed point $(\bar p,\bar q) \in \Omega$, which satisfies \eqref{eq3.2}. Thus $(\bar p,\bar q)$ is a Nash equilibrium point of the game. 
\hfill $\h$

Next we use the above variational arguments, with applying an appropriate fixed-point theorem, to establish the existence of Nash equilibrium for two-person games in Hilbert spaces.

\begin{Theorem}[\bf existence of Nash equilibrium in Hilbert spaces]\label{nash-hilb}
Let $H_1$ and $H_2$ be Hilbert spaces. Consider a two-person game $\{\Omega_{i}, u_{i}\}$ for $i=1,2$, where $\Omega_{1} \subset H_1$ and $\Omega_{2} \subset H_2$ are nonempty compact convex sets. Let the payoff functions $u_{i}\colon  \Omega_{1} \times \Omega_{2} \rightarrow \mathbb{R}$ be given by
$$
u_{1}(p, q): = \langle p,Aq \rangle\;\mbox{ and }\;u_{2}(p, q): = \langle q,Bp \rangle.
$$
where $A\colon  H_2 \to H_1$ and $B\colon  H_1 \to H_2$ are continuous linear operators. Then this game admits a Nash equilibrium.
\end{Theorem}
{\bf Proof.} Arguing as in the proof of Theorem~\ref{theo3.3} with employing the optimality condition from Proposition~\ref{prop2.4}, which holds in Hilbert spaces, brings us to the inclusions
$$	
A\bar q\in N(\bar p;\Omega_{1}) \ \; \text{and} \ \;  B\bar p\in N(\bar q;\Omega_{2}),
$$
i.e., $(A\bar q, B\bar p) \in N(\bar p;\Omega_{1}) \times N(\bar q;\Omega_{2})$. Then we proceed similarly to the proof of 
Theorem~\ref{theo3.3} with the usage of Proposition~\ref{prop2.7} held in Hilbert spaces and applying Schauder's fixed-point result from Theorem~\ref{theo2.8} instead of its Brouwer's counterpart in finite dimensions. This verifies the existence of a Nash equilibrium point.\hfill $\h$

\section{Nash Equilibrium with Nondifferentiable Payoff Functions in LCTV Spaces}

In this section, we develop our variational approach to establish the existence of Nash equilibrium points for two-person games in LCTV spaces, where partially concave payoff functions may not be differentiable. As a consequence of this main result, a general version of von Neumann's minimax theorem for concave-convex functions is obtained in the LCTV setting.

First we present a simple lemma used in the proof of the main result below.

\begin{Lemma}[\bf properties of stationary multifunctions]\label{lem4.1}
Let $X, Y$ be LCTV spaces, and let $f\colon X \times Y \to (-\infty,\infty]$ be a function such that $f(\cdot, y)$ is convex for every $y\in Y$.  Define the stationary multifunction $G\colon  Y \tto X$ associated with $f$ by
$$
G(y): = \big\{x \in X \; \big |\; 0 \in \partial_x f(x,y)\big\}.
$$
Then $G$ satisfies the following properties:

{\bf(a)} The set $G(y)$ is convex for any fixed $y \in Y$.

{\bf(b)} The graph of $G$ is closed if $K:=\dom f$ is closed and $f$ is continuous on $K$.
\end{Lemma}
{\bf Proof.} To verify (a), fix $y\in Y$ and pick any $x_1, x_2 \in G(y)$ together with $0 < \lambda < 1$. Then we have the subgradient inclusions
$$
0 \in \partial_x f(x_1,y) \ \; \mbox{\rm and }\ 0 \in \partial_x f(x_2,y).
$$
The convexity of $f(\cdot,y)$ and the definition of the partial subdifferential yield
\begin{equation*}
0\leq  f(x,y) - f(x_1,y) \ \; \mbox{\rm and }\ 0\leq  f(x,y) - f(x_2,y)\ \; \mbox{\rm for all }\ x\in X.
\end{equation*}
Multiplying the first inequality above by $\lambda$, the second by $1-\lambda$, and then summing up the resulting relationships give us
\begin{equation} \label{eq4.1}
\begin{aligned}
0   &\leq  \lambda\big(f(x,y) - f(x_1,y)\big)+ (1-\lambda)\big(f(x,y) - f(x_2,y)\big)\\
    &=f(x, y)-\big( \lambda f(x_1,y) + (1-\lambda)f(x_2,y)\big)\ \; \mbox{\rm for all }\ x\in X.
\end{aligned}
\end{equation}
Using again the convexity of $f(\cdot,y)$, we get 
\begin{equation*}
f(\lambda x_1 + (1-\lambda) x_2,y)\leq  \lambda f(x_1,y) + (1-\lambda)f(x_2,y), 
\end{equation*}
and then combining the latter with \eqref{eq4.1} tells us that
$$
0\leq  f(x,y) - f\big(\lambda x_1 + (1-\lambda) x_2,y\big)\ \; \mbox{\rm for all }\ x\in X.
$$
Therefore, $0 \in \partial_x f(\lambda x_1 + (1-\lambda) x_2,y)$, and hence $\lambda x_1 + (1-\lambda) x_2 \in G(y)$, which justifies the convexity of $G(y)$.

Next we verify (b). Take any net $\{(x_\tau,y_\tau)\}_{\tau \in T}$ such that $x_\tau \in G(y_\tau)$ for all $\tau \in T$ and $(x_\tau,y_\tau) \to (x_0,y_0)$. We aim at showing that $x_0 \in G(y_0)$. Indeed,  $x_\tau \in G(y_\tau)$ means that $0 \in \partial_x f(x_\tau,y_\tau)$, and hence $(x_\tau, y_\tau)\in \dom f$ with
$$
0 \leq  f(x,y_\tau) - f(x_\tau,y_\tau) \ \; \mbox{\rm for all } ~ x \in X.
$$
Passing to a limit above and taking into account the continuity of $f$ on the closed set $K$ leads us to $(x_0, y_0)\in \dom f$ and
$$
0 \leq  f(x,y_0) - f(x_0,y_0)\ \; \mbox{\rm whenever } ~ x \in X,
$$
which shows that $x_0 \in G(y_0)$ and thus completes the proof of the lemma.
\hfill $\h$

Now we are ready to employ our variational device to establishing the existence of Nash equilibrium for nonsmooth two-person games in LCTV spaces. To the best of our knowledge, in such generality this result has never been achieved before.

\begin{Theorem}[general Nash equilibrium theorem in LCTV spaces] \label{theo4.2}
 Let $X_1, X_2$ be LCTV spaces,  and let $u_{i}\colon X_1 \times X_2  \rightarrow \mathbb{R}$ for $i=1, 2$ be continuous payoff functions in the two-person game $\{\Omega_{i}, u_{i}\}$, where $\Omega_{1} \subset X_1$ and $\Omega_{2} \subset X_2$ are nonempty compact convex sets. Suppose that the functions $u_{i}$ for $i=1, 2$ have the following properties:

{\bf(a)} $u_1(\cdot,q)$ is concave for every fixed $q \in \Omega_2$.

{\bf(b)} $u_2(p,\cdot)$ is concave for every fixed $p \in \Omega_1$.

Then this game admits a Nash equilibrium point.
\end{Theorem}
{\bf Proof.} Remember that for $(\bar p,\bar q) \in \Omega_{1} \times \Omega_{2}$ being a Nash equilibrium point means that both conditions in \eqref{nash1} are satisfied. The optimality condition of Proposition~\ref{prop2.4} allows us to equivalently describe \eqref{nash1} via the inclusions
\begin{equation}\label{sum rule}
\begin{array}{ll}
&0 \in \partial_p (-u_1)(\bar p,\bar q) + N(\bar p;\Omega_1), \\
&0 \in \partial_q (-u_2)(\bar p,\bar q) + N(\bar q;\Omega_2).
\end{array}
\end{equation}
Considering the extended-real-valued functions $f_1, f_2: X_1 \times X_2 \to(-\infty,\infty]$ defined by $f_1 := -u_1 + \delta(\cdot,\Omega_1 \times X_2)$ and $f_2 := -u_2 + \delta(\cdot,X_1 \times \Omega_2)$, and employing the subdifferential sum rule (Moreau-Rockafellar theorem) of convex analysis, we equivalently rewrite \eqref{sum rule} as
\begin{equation} \label{eq4.3}
0 \in \partial_p f_1(\bar p,\bar q) \ \; \text{and} \ \; 0 \in \partial_p f_2(\bar p,\bar q).
\end{equation}
Define further the set-valued mappings associated with each player by
\begin{align*}
&F_1(q) := \{p \in \Omega_1 \mid 0 \in \partial_p f_1(p,q)\} \ \; \text{for } ~ q \in \Omega_2, \\
&F_2(p) := \{q \in \Omega_2 \mid 0 \in \partial_p f_2(p,q)\} \ \; \text{for } ~ p \in \Omega_1
\end{align*}
and consider the set-valued mapping $F\colon \Omega_1 \times \Omega_2 \tto \Omega_1 \times \Omega_2$ given by $F(p,q) = F_1(q) \times F_2(p)$. Observe that if $F$ has a fixed point $(\bar p,\bar q) \in F(\bar p,\bar q)$, then \eqref{eq4.3} holds, and hence $(\bar p,\bar q)$ verifies the existence of Nash equilibrium. We intend to apply Kakutani's fixed-point result of Theorem~\ref{theo2.9} to ensure that the multifunction $F$ admits a fixed point.

To furnish this, fix $q \in \Omega_2$ and deduce from the continuity of $u_1(\cdot, q)$ that this function achieves its maximum on the compact set $\Omega_1$ at some point $\bar p\in\Omega_1$. It follows from the structure of $F_1$ and the first-order optimality condition in Proposition~\ref{prop2.4} that $\bar p\in F_1(q)$, and thus the set $F_1(q)$ in nonempty. Note that $\dom f_1 = \Omega_1 \times X_2$ is a closed set. Applying Lemma~\ref{lem4.1} with $f = f_1$, we deduce that $F_1(q) = G(q) \cap \Omega_1$ in the notation of the lemma. The concavity of $u_1(\cdot, q)$ and the convexity of $\Omega_1$ yield the convexity of $f_1(\cdot, q)$ and hence the convexity of $F_1(q)$. Furthermore, the continuity of $u_1$ and the closedness of $\Omega_1$ imply that $f_1$ is continuous on $\Omega_1\times X_2$. Therefore, the set $\gph G$ is closed and so is $\gph F_1 = \gph G \cap (X_2\times \Omega_1)$. 

The case of $F_2$ is completely symmetric to the above consideration for $F_1$. Since both components of the Cartesian product $F=F_1\times F_2$ have the closed graphs and nonempty convex values, $F$ also satisfies these properties. Therefore, all the requirements of Theorem~\ref{theo2.9} are fulfilled, and the claimed existence of a Nash equilibrium point has been justified.
\hfill $\h$

As a consequence of the existence of Nash equilibrium established in Theorem~\ref{theo4.2}, we now derive a generalized version of von Neumann's minimax theorem in LCTV spaces; see further extensions via different approaches in the next section.

\begin{Corollary}[minimax equality for convex-concave functions] \label{theo4.3}
Let $X_1, X_2$ be  LCTV spaces, and let $\Omega_1 \subset X_1$,  $\Omega_2 \subset X_2$ be nonempty compact convex sets. Assume that $f\colon X_1 \times X_2 \to \mathbb{R}$ is a continuous function that is concave-convex, i.e.,

{\bf(a)} $f(\cdot,y)\colon X_1 \to \mathbb{R}$ is concave for every fixed $y \in \Omega_2$,

{\bf(b)}  $f(x,\cdot)\colon X_2 \to \mathbb{R}$ is convex for every fixed $x \in \Omega_1$.

Then we have the minimax equality
\begin{equation}\label{eq4.4}
\max_{x \in \Omega_1}\min_{y \in \Omega_2} f(x,y) = \min_{y \in \Omega_2}\max_{x \in \Omega_1} f(x,y).
\end{equation}
\end{Corollary}
{\bf Proof.}
First note that the continuity of $f$ on the compact set $\Omega_1 \times \Omega_2$ ensures that maximin and minimax in \eqref{eq4.4} are attainable. For any $x_0 \in \Omega_1$ and $y_0 \in \Omega_2$, we always have the conditions
$$
\min_{y \in \Omega_2} f(x_0, y) \leq f(x_0,y_0) \leq \max_{x \in \Omega_1} f(x,y_0),
$$
which yield therefore the inequalities
\begin{equation} \label{eq4.5}
\max_{x \in \Omega_1} \min_{y \in \Omega_2} f(x, y) \leq f(x_0,y_0) \leq \min_{y \in \Omega_2} \max_{x \in \Omega_1} f(x,y).
\end{equation}
To verify the reverse inequalities in \eqref{eq4.5}, define $u_1: = f$ and $u_2: = -f$. It follows from the assumptions on $f$ that both conditions (a) and (b) imposed in Theorem~\ref{theo4.2} are satisfied. Consequently, a Nash equilibrium pair $(\bar x,\bar y) \in \Omega_1 \times \Omega_2$ exists meaning the fulfillment of \eqref{nash1}, which is equivalent to
\begin{align*}
& f(x,\bar y) \leq f(\bar x,\bar y) \ \; \text{for all} ~ x \in \Omega_1, \\
& f(\bar x, y) \geq f(\bar x,\bar y) \ \; \text{for all } ~ y \in \Omega_2.
\end{align*}
Therefore, we have the relationships
$$
\max_{x \in \Omega_1} f(x,\bar y) \leq f(\bar x,\bar y) \leq \min_{y \in \Omega_2} f(\bar x,y),
$$
from which it follows that
\begin{equation} \label{eq4.6}
\min_{y \in \Omega_2} \max_{x \in \Omega_1} f(x,y) \leq f(x_0, y_0) \leq \max_{x \in \Omega_1} \min_{y \in \Omega_2} f(x, y).
\end{equation}
Combining finally the inequalities in \eqref{eq4.5} and \eqref{eq4.6}, we arrive at the minimax equality \eqref{eq4.4} and thus complete the proof of the corollary.
\hfill $\h$

\section{Extended Minimax Theorems in LCTV Spaces}

In this section, we revisit minimax theorems and obtain extensions of Corollary~\ref{theo4.3} by using tools of convex analysis and fixed-point theory in LCTV spaces. Proceeding in this way allows us to relax continuity and convexity-concavity assumptions imposed in Corollary~\ref{theo4.3}. We refer the reader to the books  \cite{borwein05,simons,zalinescu02} and the bibliographies therein for various results in this direction given in diverse infinite-dimensional settings, where the proofs are  different from those developed below.  

Recall that the (unilateral) Weierstrass theorem guarantees that any upper semicontinuous function attains its global maximum on any nonempty compact subset $\Omega$ of an LCTV space $X$, i.e., there exists $x_0 \in \Omega$ such that
$$
f(x) \leq f(x_0) \ \; \text{for all } ~ x \in \Omega.
$$

Here is the first minimax theorem of this section, which does not impose the continuity assumption on the function under consideration, while does not assert that the minimum  with respect to a compact set is attained.  

\begin{Theorem}[\bf inf-max relationship under upper semicontinuity]\label{theo5.1} Let $X_1, X_2$ be LCTV spaces, and let $\Omega_i \subset X_i$, $i=1,2$, be nonempty compact convex sets. Consider a function $f\colon\Omega_1 \times \Omega_2 \to \mathbb{R}$ satisfying the following conditions:

{\bf(a)} $f(\cdot,y)$ is concave and upper semicontinuous for every fixed $y \in \Omega_2$.

{\bf(b)} $f(x,\cdot)$ is convex for every fixed $x\in\Omega_1$.

Then we have the inf-max equality
\begin{equation} \label{eq5.1}
 \max_{x \in \Omega_1}\inf_{y \in \Omega_2} f(x,y) = \inf_{y \in \Omega_2}\max_{x \in \Omega_1} f(x,y).
\end{equation}
\end{Theorem}
{\bf Proof.} The maximum value on the right-hand side of \eqref{eq5.1} is attained by the Weierstrass existence theorem. To verify that the maximum is achieved also on the left-hand side therein, we are going to show that the value function
$$
h(x):= \inf_{y \in \Omega_2} f(x,y),\ \;  x\in \Omega_1,
$$
is upper semicontinuous on $\Omega_1$. To furnish this, for any real number $\alpha$ consider the upper level set $U_{\alpha}(h) = \{x \in \Omega_1 \mid h(x) \geq \alpha\}$, which can be expressed as
\begin{align*}
U_{\alpha}(h)
&= \big\{x \in \Omega_1 \; \big |\; f(x,y) \geq \alpha\;\mbox{ for all }\;y \in \Omega_2\big\} \\
&= \bigcap_{y \in \Omega_2} \big\{x \in \Omega_1 \; \big | \; f(x,y) \geq \alpha\big\}.
\end{align*}
The upper semicontinuity of $f(\cdot,y)$ ensures that each set $\{x \in \Omega_1 \mid f(x,y) \geq \alpha\}$ is closed. Since the intersection of closed sets is closed, the upper level set $U_{\alpha}(h)$ is proved to be closed. Therefore, the value function $h$ is upper semicontinuous on $\Omega_1$.

It remains to verify the inequality ``$\ge$" in \eqref{eq5.1} because the reverse one is obvious. Denote $\sigma:=\max_{x\in\Omega_1}\inf_{y\in\Omega_2}f(x,y)$ and observe that $\sigma\in[-\infty,\infty)$. For any real number $\alpha>\sigma$, we 
get the condition
$$
\inf_{y \in \Omega_2} f(x, y) < \alpha\;\mbox{ for all }\;x\in\Omega_1.
$$
Given further any $y \in \Omega_2$, define the sets $A(y): = \{x \in \Omega_1 \mid f(x,y) \geq \alpha\}$ and observe that each set $A(y)$ is closed by the upper semicontinuity of $f(\cdot,y)$. Moreover, $\bigcap_{y \in \Omega_2} A(y) = \emptyset$. The compactness of $\Omega_1$ allows us to find finitely many elements $y_1,\ldots,y_n \in \Omega_2$ such that $\bigcap_{j = 1}^n A(y_j) = \emptyset$. In other words,
\begin{equation} \label{eq5.2}
\min_{j\in\{1,\ldots,n\}} f(x, y_j) < \alpha\;\mbox{ for every }\;x\in\Omega_1.
\end{equation}
Define now the mapping $\varphi\colon \Omega_1 \to \mathbb{R}^n$ by $\varphi(x) = (\varphi_1(x),\ldots,\varphi_n(x))$, where
$$
\varphi_j(x): = f(x, y_j) - \alpha \; \ \text{ whenever } \; \ j = 1,\ldots,n.
$$
We claim that $\conv(\varphi(\Omega_1)) \cap \mathbb{R}^n_+ = \emptyset$. To verify this, take any $\lambda_1, \ldots, \lambda_m \in [0, 1]$ with $\sum_{i=1}^m {\lambda_i = 1}$ and $x_1, \ldots, x_m \in \Omega_1$ for $m\in \N$. The concavity of $f(\cdot, y)$ yields
$$
f\Big(\sum_{i=1}^m \lambda_ix_i, y_j\Big) \geq \sum_{i=1}^m  \lambda_if(x_i,y_j)
$$
for all $j=1, \ldots, n$. This implies therefore that
$$
\varphi_j\Big(\sum_{i=1}^m \lambda_ix_i\Big) \geq \sum_{i=1}^m  \lambda_i\varphi_j(x_i).
$$
Since $\Omega_1$ is convex, we see that $\sum_{i=1}^m \lambda_ix_i \in \Omega_1$. Condition \eqref{eq5.2} ensures the existence of $j_0$ such that $\varphi_{j_0}(\sum_{i=1}^m \lambda_ix_i) < 0$, which implies in turn that $\sum_{i=1}^m  \lambda_i\varphi_{j_0}(x_i) < 0$. Hence, $\sum_{i=1}^m  \lambda_i\varphi(x_i)\notin \mathbb{R}^n_+$. By the convex separation theorem, there exists a vector $(\beta_1, \ldots, \beta_n) \in \mathbb{R}^n \setminus \{0\}$ for which
\begin{align}
&\sum_{j=1}^n \beta_j \varphi_j(x) \leq 0\;\mbox{ whenever }\;x \in \Omega_1, \label{eq5.3} \\
&\sum_{j=1}^n \beta_j e_j \geq 0\;\mbox{ whenever }\; (e_1, \ldots, e_n) \in \mathbb{R}^n_+.\label{eq5.4}
\end{align}
It follows from inequality \eqref{eq5.4} that $\beta_j \geq 0$ for all $j$. Dividing both sides of \eqref{eq5.3} by $\sum_{j=1}^n \beta_j$, we may assume that $\sum_{j=1}^n \beta_j = 1$. Consequently, \eqref{eq5.3} yields
$$
f\Big(x, \sum_{j=1}^n \beta_j y_j\Big) \leq \sum_{j=1}^n \beta_j f(x, y_j) \leq \alpha\;\mbox{ for all }\; x \in \Omega_1,
 $$
where the first inequality results from the convexity of $f(x, \cdot)$. Since $\Omega_2$ contains the sum $\sum_{j=1}^n \beta_j y_j$ due to the set convexity, we get
$$
\inf_{y \in \Omega_2} \max_{x \in \Omega_1} f(x, y) \leq \max_{x \in \Omega_1} f\Big(x, \sum_{j=1}^n \beta_j y_j\Big) \leq \alpha.
$$
Letting $\alpha \to \sigma$ tells us that
$$
\inf_{y \in \Omega_2} \max_{x \in \Omega_1} f(x, y) \leq \sigma = \max_{x \in \Omega_1} \inf_{y \in \Omega_2} f(x, y),
$$
which thus completes the proof of the theorem.
\hfill $\h$

If in addition we impose the lower semicontinuity of $f$ with respect to the first variable, then we get a direct extension of the minimax result of Corollary~\ref{theo4.3}.

\begin{Corollary}[\bf minimax equality without continuity] \label{cor5.2} Let $X_i$ and $\Omega_i$, $i=1,2$, be as in Theorem~{\rm\ref{theo5.1}}, and let $f\colon  \Omega_1 \times \Omega_2 \to \mathbb{R}$ satisfy the following conditions:

{\bf(a)} $f(\cdot,y)$ is concave and upper semicontinuous for every fixed $y \in \Omega_2$.
 
{\bf(b)} $f(x,\cdot)$ is convex and lower semicontinuous for every fixed $x \in \Omega_1$.

Then we have the minimax equality \eqref{eq4.4}. Furthermore, the function $f$ admits a saddle point $(\bar x, \bar y) \in \Omega_1 \times \Omega_2$.
\end{Corollary}

{\bf Proof.} It follows from the imposed lower semicontinuity that the infima in \eqref{eq5.1} become minima. Consequently, there exists a pair \((\bar{x}, \bar{y}) \in \Omega_1 \times \Omega_2\) such that
$$
\min_{y \in \Omega_2} f(\bar{x}, y) = \max_{x \in \Omega_1} f(x, \bar{y}).
$$
Due to the inequalities
$$
\min_{y \in \Omega_2} f(\bar{x}, y) \leq f(\bar{x}, \bar{y}) \quad \text{and} \quad \max_{x \in \Omega_1} f(x, \bar{y}) \geq f(\bar{x}, \bar{y}),
$$
this equality can be expressed as:
$$
\min_{y \in \Omega_2} f(\bar{x}, y) = f(\bar{x}, \bar{y}) = \max_{x \in \Omega_1} f(x, \bar{y}).
$$
Therefore, we arrive at the conditions
$$
f(\bar x,y) \geq f(\bar x,\bar y) \geq f(x,\bar y)\;\mbox{ for all }\;(x,y) \in \Omega_1 \times \Omega_2
$$
meaning that $(\bar{x}, \bar{y})$ is a saddle point of the game.
\hfill $\h$

\section{Conclusions}\label{sec:conc}

In this paper, we develop a new variational approach of convex analysis and generalized differentiation, which brings us to establish the existence of Nash equilibrium in a general setting of LCTV spaces. In this way, we also derive various generalized versions of von Neumann's minimax theorem in LCTV spaces without continuity and convexity-concavity assumptions on the function in question. 

Among unsolved issues for our future research, the two major topics should be mentioned. First we plan to investigate the possibility of extending Theorem \ref{theo4.2} on the existence of Nash equilibrium beyond the continuity and convexity-concavity assumptions on payoff functions in zero-sum games. The second topic revolves around further extensions to non-zero sum games where, in particular, the existence of saddle points may not imply the existence of Nash equilibrium.

\end{document}